 \documentclass[12pt,leqno]{amsart}
\newtheorem{dummy}{}[section]

\newtheorem{question}[dummy]{Question}
\begin{document}
\bibliographystyle{plain}
\title{Derived Langlands V:The simplicial and Hopflike categories } 
\author{Victor P. Snaith}
\date{30 September 2020}
\maketitle
 \tableofcontents 
 
 \section{Introduction}
 
 This is going to be a rather strange\footnote{The strangeness is a consequence of a medical time-table beyond my control \cite{JMOS18}. I shall try my best to assemble enough of a coherent bunch of slightly non-classical Langlands notions, enough perhaps to interest the occasional reader.} introduction in which I am going to attempt to synopsise the rather sketchy synthesis which is presented in my Derived Langlands series consisting of a research monograph and four essays (\cite{Sn18}, \cite{Sn20}, \cite{Sn20b}, \cite{newindnotes}, \cite{Sn20c}).
 
 In the 1980's homotopy theorists and algebraic topologists became interested in the Burnside ring of finite groups $G$ (and later $G$ became compact Lie) in connection with the Segal conjecture. Denotes by 
 $\Omega(G)$ the Burnside ring is the ring of finite $G$ subsets or of permutation representations. The interest was  a theorem which described the $\Omega(G)$-adic completion in terms of stable cohomotopy groups of the classifying space $BG$ of $G$.
 
 Around 1985 I became aware of a problem (inspired by Dwork's work on the local function equation and local root numbers) which was occupying Jean-Pierre Serre (footnote p.71 \cite{Ser77}), Bob Langlands and Pierre Deligne. The problem was to find the presentation for the additive complex representation of $G$ in terms of generators given by induced representation of the form ${\rm Ind}_{H}^{G}(\phi)$ where $\phi$ is a one-dimensional representation of $H$.
 
 Analogous to $\Omega(G)$ one can define $R_{+}(G)$ \cite{Sn94} as a ring generated by conjugacy classes of $(H, \phi)$'s. In 1985 I learned how to identify the $\Omega(G)$-adic completion of $R_{+}(G)$ in terms of stable cohomotopy \cite{Sn89}. There is also a completion theorem concerning $R(G)$, due to Atiyah. The answer to the question of (footnote p.71 \cite{Ser77}) was easy to provide after these completions because there is a stable homotopy transfer from Atiyah's completion to that of $R_{+}(G)$ which splits the induction map in the opposite direct, which is surjective by Brauer's induction theorem. To remove the completion, once one knows the formula after completion, given in \cite{Sn89}, can be accoplished in several ways - I used the Atiyah-Singer Index Theorem in \cite{Sn88}.
 
 Now I can begin to describe what I believe to be the germs of a coherent synthesis in which I introduce a number of analogues as we expand from the context of finite dimensional representations of finite groups to various notions of ``admissible'' representations of locally $p$-adic Lie groups. 
 
 Although ad hoc chain complexes, derived categories and the notion of a ``monomial resolution'' of a finite dimensional representation of a finite group were evident through the algebraic topology of \cite{Sn88}
 the proper categorical setting was only finally described in \cite{Bol01}. The construction is performed in the manner of the Freyd-Mitchell Theorem which passes from an additive category, like the monomial category for finite groups, to a functor category which is an abelian module category over a ring of monomial endomorphisms given by generators and relation which appear, extended to the $p$-adic Lie group context, in 
\cite{Sn18} - in a clumsy manner - and in \cite{Sn20}. I call this ring the hyperHecke algebra of $G$.

Even in the finite group case we homotopy theorists had overlooked the hyperHecke  and its action on $R_{+}(G)$ even though we were very enthusiastic about completion theorems which turned the latter into stable cohomotopy! 

The problem with the Langlands Programme was that there was an immense amount of technical background, familiar to the experts from whom I was isolated once I moved back from Canada to the UK. However, in 2005, it dawned on me that the monomial resolution of \cite{Bol01} was - on the nose via the  Freyd-Mitchell  procedure - the bar resolution of the hyperHecke algebra. 

Now for a bit more of the synthesis. When $G$ moves on to being a locally $p$-adic Lie group then $R{+}(G)$, by the choice made in \cite{Sn18} and \cite{Sn20}, becomes a ring of monomial isomorphism classes of
pairs $(H, \phi)$ where $H$ is compact open mod centre, containing the centre of $G$, and $\phi$ is a continuous character into some algebraic closed field. In the classical Langlands situation the field has a different residue character from $p$ and these characters have finite image when restricted to compact, open subgroups. The hyperHecke algebra
\cite{Sn20} is generated by triples $((K, \psi), g , (H, \phi))$ modulo relations (\cite{Sn18}, \cite{Sn20}).

The monomial resolution of an admissible representation $V$ - (sometimes I tended to concentrate the case where $V$ had a fixed central character $\underline{\phi}$ which coincided with the restrictions of the $\psi$'s and $\phi$'s to the centre of $G$ - is then given as follows. Take Tammo tom Dieck's space for $G$ and the family of $(H, \phi)$'s as above - see final Appendix of \cite{Sn18}. In \cite{Sn18} I overcomplicated things until I finally realised that, as far as equivariant topology was concerned, Tammo's space generalise the Bruhat-Tits buildings. It is an equivariant cell complex in which the stabliliser of any simplex is one of our $H$'s. Given admissible $V$ and a simplex in Tammo's space the bar resolution of the restriction of $V$ to the simplex stabiliser gives a local coefficient system (i.e. a sheaf) whose associated double complex - through the 
 Freyd-Mitchell  routine - is, according to \cite{Sn20}, a monomial complex in the derived category of the monomial category of $G$.
 In other words an embedding of admissible representation objects in the Langlands Programme into a derived category of an (slightly awkwardly non-abelian) additive category.
 
 Next recall from \cite{PD84} (see also \cite{Sn20} \S13) how idempotented algebras and convolution algebras enter into the classical Langlands Programme. In particular, in order to describe the Bernstein centre of the category of classical admissible representations of $G$ one identifies the representation category with a category of Hecke modules over the convolution algebra given by the Hecke algebra. In order to substantiate the synthesis I gave a formula in (\cite{Sn20} \S7) for the product in the hyperHecke algebra at least in the classical case when characters on compact open subgroups have finite image. I believe that the same formula will make sense with the aid of some version of Mahler's Theorem \cite{AMR2000}. The point about the hyperHecke algebra is that it is one of the analogous notions which appears once one takes seriously not merely the subgroups $H$ but also the characters $\phi$ on them. The classical Hecke algebra is analogous to the sub-convolution algebra in which the $\phi$'s are trivial. Relating the product in the hyperHecke algebra to convolution products suggests, I believe, that the hyperHecke algebra is a generalisation in the right direction.
 
 In addition, in (\cite{Sn20} \S3) I was able to give a combinatorial description of the centre of the hyperHecke algebra and of various centralising subgroups, with a view to developing a hyperHecke algebra approach to the Bernstein centre of the representation category via monomial resolutions. The idea being that subrings of the hyperHecke algebra generated by partially central elements induce ``central endomorphisms'' of the monomial resolution which, in turn, would induce morphisms in the Bernstein centre of the representation category.
 
 In (\cite{Sn20} \S3) I define the notion of ${\mathcal M}_{cmc, \underline{\phi}}(G)$-admissible representations. In the classical Langlands situation this coincides with the usual definition 
 (\cite{Sn20} Proposition 4.1). However, for a locally $p$-adic Lie group and representations over an algebraically closed local field of residue characteristic $p$ (the di-$p$-adic situation) this is a different concept. 

This type of admissible representation has an appropriate version of induction - called 
${\mathcal M}_{c}(G)$-smooth induction - which is defined in (\cite{newindnotes} Definition 1.2; \S1.5 for the $c$-IND version). This induction possesses analogues of the usual short exact sequence and categorical properties (\cite{newindnotes} \S\S1.6-1.13).
 
 In (\cite{newindnotes}  \S2 Theorem 2.2) I claim, for this type of admissibility and  ${\mathcal M}_{c}(G)$-smooth induction, to establish the analogue of Jacquet's Theorem 
 (\cite{AJS79} Chapter 2, \S2.3) for ${\mathcal M}_{cmc, \underline{\phi}}(G)$-admissible representations.
 In view of the fundamental role of Jacquet's Theorem in the classical classification of admissible representations, if correct, this property should have classificatory consequences for 
 ${\mathcal M}_{cmc, \underline{\phi}}(G)$-admissibility.
 
 Finally we come to Hopf algebras, PSH algebras and Hopflike categories, which entered into the classification of classical admissible representations in \cite{BZ76}, \cite{BZ77} 
(and \cite{AVZ81} for finite general linear groups). For example, in the finite case the integral Milnor-Moore structure theorem, for the PSH algebra $\oplus_{n} \ R(GL_{n}{\mathbb F}_{q})$, simply gives all the complex irreducibles in terms of the primitives (aka cuspidals) in the PSH algebra. The locally $p$-adic $GL_{n}$ is considerably more technical but the basic idea is similar (\cite{BZ76}, \cite{BZ77}).

Now, even in the finite case,  $\oplus_{n} \ R_{+}(GL_{n}{\mathbb F}_{q})$ is not a PSH algebra because 
$R_{+}(GL_{n}{\mathbb F}_{q} \times GL_{m}{\mathbb F}_{q} )$ is not isomorphic to 
$R_{+}(GL_{n}{\mathbb F}_{q}) \otimes R_{+}GL_{m}{\mathbb F}_{q} )$. However, my idea is that perhaps
\[   \begin{array}{l}
 \oplus_{(a_{1}, a_{2}, \ldots , a_{k}) } \ R_{+}(  GL_{a_{1}}{\mathbb F}_{q}  \times 
 GL_{a_{2} }{\mathbb F}_{q} 
 \times \ldots \times GL_{a_{k}}{\mathbb F}_{q}  ) ,     \\
 \end{array}  \]
 indexed over ordered, non-negative integer partitions\footnote{Often  called compositions in the combinatorial literature, apparently.} might be a PSH algebra - and for combinatorial reasons which apply equally to the locally $p$-adic case. Were this true, it might  even be true also for a similar direct sum of hyperHecke algebras of products of $p$-adic $GL_{n}$'s.
 
 In this essay, I can only point to some incomplete evidence for these hopes and to some potential consequences. This is especially true in \S9 (``Low dimensional etc'').

We have the usual multiplication
\[ (R(a_{1}) \otimes \ldots \ \otimes R(a_{k})) \otimes   (R(b_{1}) \otimes \ldots \ \otimes R(b_{k}))
\longrightarrow  (R(a_{1} + b_{1}) \otimes \ldots \ \otimes R(a_{k} + b_{k}))  \]
and comultiplication
\[ \begin{array}{c}
 (R(a_{1}) \otimes \ldots \ \otimes R(a_{k}))  \\
\\
m^{*}  \downarrow  \\
\\
\Sigma \Sigma\ldots \Sigma  \ (R(a_{1}- e_{1}) \otimes \ldots \ \otimes R(a_{k}- e_{k})) \otimes (R(e_{1}) \otimes \ldots \ \otimes R(e_{k}))
\end{array} \]
which satisfies the Hopf condition.

Let us say that the bidegree of $R(a) \otimes (R(b) \otimes R(c))$ is $(1,2)$ and that the bidegree of 
\[ ( R(a_{1}) \otimes R(a_{2})  \otimes ((R(b) \otimes R(c)) \otimes R(d) \otimes Re)) \]
is equal to $(2,4)$ and so on. In \S9  (``Low dimensional etc'') I have define $m$ in bidegree $(1,2)$ and I believe I show how to make a modification to $m^{*}$ to give a comultiplication in bidegree $(1,2)$ (and by ``symmetry'' $(2,1)$) which satisfies the Hopf condition
\[   (m \otimes  m)({\rm shuffle})(m^{*} \otimes m^{*} \otimes m^{*} = m^{*} m . \]
The modification is entirely combinatorial and is determined by the Hopf condition. The Hopf condition will become very complicated, but let us sppose pro tem that the process succeeds in making a Hopf algebra - even a PSH (positive selfadjoint Hopf algebra) - then we are in the classification business in general - including the di-$p$-adic settings. If my wishful thinking gets us this far one might hope that a similar massive diret sum of hyperHecke algebras of products of $p$-adic general linear groups, indexed by ordered
partitions, also becomes a Hopflike structure (Hopf algebra or PSH) by means of the candidate structure whih I 
give in (\cite{Sn20b} \S12).

 If this sort of Hopf algebra works it should dove-tail nicely with the Langlands `unique irreducible quotient" \cite{AJS78} result and maybe even the modular case of Florian Herzig \cite{FH2011 } and and its sequel with collaborators\footnote{There is a slightly relevant anecdotal topological footnote on Question \ref{8.1}}.

 \section{Simplicial sets}

 Let $\Delta$ denote the simplex category. The objects of $\Delta$ are nonempty linearly ordered sets of the form
\[  [n] = {0, 1, ..., n} \]
with $n \geq 0$. The morphisms in $\Delta$ are (non-strictly) order-preserving functions between these sets.
A simplicial set $X$ is a contravariant functor
\[  X : \Delta \longrightarrow   {\rm Set} \]
where ${\rm Set}$ is the category of sets. (Alternatively and equivalently, one may define simplicial sets as covariant functors from the opposite category $\Delta^{op}$ to ${\rm Set}$.) Simplicial sets are therefore nothing but presheaves on $\Delta$. Given a simplicial set $X$, we often write $X_{n}$ instead of $X([n])$.

Simplicial sets form a category, usually denoted ${\rm sSet}$, whose objects are simplicial sets and whose morphisms are natural transformations between them.

If we consider covariant functors $X : \Delta \longrightarrow {\rm Set}$ instead of contravariant ones, we arrive at the definition of a cosimplicial set. The corresponding category of cosimplicial sets is denoted by ${\rm cSet}$.

\underline{{\bf Face and degeneracy maps}}

The simplex category $\Delta$ is generated by two particularly important families of morphisms (maps), whose images under a given simplicial set functor are called face maps and degeneracy maps of that simplicial set.

The face maps of a simplicial set $X$ are the images in that simplicial set of the morphisms 
\[  \partial^{n,0} , \partial^{n,1}  , \ldots , \partial^{n, n} : [ n -1 ] \longrightarrow [ n ]  \]
where $\partial^{n,i} $ is the only (order-preserving) injection $[ n - 1 ] \longrightarrow [ n ]$
which "misses" $i$. Let us denote these face maps by 
\[ d_{n , 0} , d_{n , 1}, \ldots , d_{ n , n}  :   X_{n}  \longrightarrow X_{ n - 1} . \]
 If the first index is clear, we just write $d_{ i}$ instead of $d_{ n , i}$. 

The degeneracy maps of the simplicial set $X$ are the images in that simplicial set of the morphisms 
\[  \sigma^{ n , 0} ,   \sigma^{ n , 1}, \ldots ,   \sigma^{ n , n} : [ n + 1 ] \longrightarrow  [ n ]  \]
 where  $ \sigma^{ n , i}$ is the only (order-preserving) surjection $[ n + 1 ] \longrightarrow [ n ]$
which "hits" $i$ twice. Let us denote these degeneracy maps by 
\[ s_{ n , 0} , s_{ n , 1}, \ldots , s_{ n , n} : X_{n} \longrightarrow X_{n+1} \]
 If the first index is clear, we just write $s_{ i}$ instead of $s_{ n , i}$.

The defined maps satisfy the following simplicial identities:
\[ \begin{array}{l}
  d_{i}d _{j} = d_{ j - 1 }d_{ i}  \  {\rm if} \  i < j,  ({\rm i.e.} \  d_{ n - 1 , i} d_{ n , j} = d_{ n - 1 , j - 1 }d_{ n , i} \
  {\rm  if } \ 0 \leq  i < j \leq  n.) \\
  \\
    d_{ i} s_{ j} = s_{ j - 1 }d_{ i} \ {\rm  if} \   i < j,   \\
    \\
    d_{ i} s_{ j} = 1 \  {\rm  if} \   i = j \ {\rm or} \  i = j + 1 , \\
    \\
    d_{ i} s_{ j} = s_{ j} d_{ i - 1} \ {\rm  if} \  i > j + 1  , \\
    \\
    s_{ i} s_{ j} = s_{ j + 1} s_{ i} \ {\rm  if} \  i \leq  j.
 \end{array} \]

Conversely, given a sequence of sets $X_{n}$ together with maps 
\[  d_{ n, i} : X_{n} \longrightarrow X_{ n -1} \
  {\rm and } \    s_{ n , i} : X_{ n} \longrightarrow  X_{ n + 1}  \]
   that satisfy the simplicial identities, there is a unique simplicial set X that has these face and degeneracy maps. So the identities provide an alternative way to define simplicial sets. 
 
 \section{Hopflike sets}
 
 Now I intend to write out the analogue of simplicial identities for the category whose objects are ordered partitions $\underline{n} = ( n_{1}, \ldots , n_{t})$ of strictly positiive\footnote{I think I shall have to allow in some zeroes and then consider two partions as equivalent if the ordered string of strictly positive integers in them are equal as ordered partitions. This at first sight appears to have added some basic morphisms to the Hopflike category but it has only added identity morphisms which were there already in the category.} integers which partition $r \geq 1$ (i.e. $r = n_{1} + \ldots + n_{t}$)\footnote{Francis Clarke helped me out with a number of typographic and mathematical errors in this section. He also pointed out that a ``modern combinatorialist''  would call these ordered partitions ``compositions'' (see https://en.wikipedia.org/wiki/Composition$\underline{\hspace{3pt}}$ (combinatorics)).}. We shall define a morphism for $i = 1,2, \ldots , t-1$
 \[   \partial^{t,i} :    ( n_{1}, \ldots , n_{t}) \longrightarrow  ( n_{1}, \ldots , n_{i-1}, n_{i} + n_{i+1}, n_{i+2}, \ldots n_{t}) \]
 and for $2 \leq n_{i}, 0 < a < n_{i}, i = 1, \ldots , t$ 
 \[ \sigma^{t,i,a} :   ( n_{1}, \ldots , n_{t}) \longrightarrow   ( n_{1}, \ldots , n_{i-1}, a , n_{i} -a,  n_{i+1}, \ldots n_{t}) \]
 
 We have a further family of morphisms between two ordered partitions of $n$  (c.f. \cite{AVZ81} p. 170) 
\[  \alpha = (a_{1}, a_{2}, \ldots , a_{r} ) \ {\rm and}  \  \beta = (b_{1}, b_{2}, \ldots , b_{s})  . \]
Suppose that we have an $r \times s$ matrix of (strictly\footnote{Not strictly, see the footnote about zeroes. Also the block sum of two matrices needs zeroes to be allowable.}) positive integers 
\[ K = \left( \begin{array}{ccccc}
k_{1,1} & k_{1,2} & \ldots & \ldots & k_{1,s}  \\
\\
k_{2,1}&k_{2,2}& \ldots & \ldots & k_{2,s} \\
\\
\vdots & \vdots &  \vdots &  \vdots &  \vdots \\
\\
k_{r,1} & k_{r,2} & \ldots & \ldots & k_{r,s} \\
\end{array} \right)  \]
such that $ \sum_{j} \ k_{i,j} = a_{i}$ for $1 \leq i \leq r$ and  $ \sum_{i} \ k_{i,j} = b_{j}$ for $1 \leq j \leq s$.
Therefore we have two ordered partitions of $n$
\[ \kappa(\alpha) :  (k_{1,1}, k_{1,2}, \ldots , k_{1,s}, k_{2,1}, \ldots , k_{2,s}, \ldots , k_{r,1}, \ldots , k_{r,s}  ) \]
and
\[ \kappa(\beta) :( k_{1,1}, k_{2,1}, \ldots , k_{r,1}, k_{1,2}, \ldots , k_{r,2}, \ldots , k_{1,s}, \ldots , k_{r,s}   ) . \]
Then we have a morphism
\[ \tau(K) : \kappa(\alpha) \longrightarrow \kappa(\beta) . \]
We envision $\tau(K)$ as the following permutation of the ordered $n$-tuple
\linebreak
 $(1,2,3, \ldots , n)$.
Think of a partition as subdividing the ordered set 
\linebreak
$(1,2,3, \ldots , n)$ into abutting intervals (called ``blocks'' in \cite{AVZ81}) of which, for $\alpha$ say, the first is $I_{1} = (1,2, \ldots , a_{1})$, the second is $I_{2} = (a_{1}+1, a_{1}+2, \ldots , a_{1}+a_{2})$ and so on up to $I_{r}$. Similarly the intervals for 
$\beta$ are $J_{1}, \ldots , J_{s}$. Then in the two partitions, $\kappa(\alpha)$ and $ \kappa(\beta)$, involving the $k_{i,j}$'s given above
there is a permutation $\sigma_{K} \in \Sigma_{n}$ which shifts (by a translation) the $k_{i,j}$-interval
of the $\kappa(\alpha)$ into the $i$-th interval of $J_{j}$ in $\kappa(\beta)$. For example, the integers in the $k_{1,2}$-interval in the upper partition goes, by a linear shift, to the interval of integers $(b_{1}+1, b_{1}+2, \ldots , b_{1}+k_{1,2})$ in the lower partition.

Next we must describe the compositional relations (analogues of the simplicial identities of \S1) between these morphisms.

Firstly the relations between compositions of $\tau(K)$'s are precisely all those relations between permutations of this type which make sense (are closed under composition).

The relations between $\partial^{t,i}$'s:
\[ \begin{array}{l}
   \partial^{t-1,j} \cdot \partial^{t,i}  =   \partial^{t-1, i-1} \cdot       \partial^{t,j}  \   {\rm for} \ j \leq i-2  \\
   \\
  \partial^{t-1,i-1} \cdot \partial^{t,i} =   \partial^{t-1,i-1 } \cdot   \partial^{t,i-1} \  {\rm for} \  j = i-1 \\
  \\
      \partial^{t-1,i} \cdot \partial^{t,i} =            \partial^{t-1,i}     \cdot   \partial^{t, i+1}               \  {\rm for} \  j = i \\
     \\
   \partial^{t-1, j-1}    \cdot \partial^{t,i}    =    \partial^{t-1,i} \cdot \partial^{t,j}                          \  {\rm for} \  j \geq i+1 \\
\end{array} \]

The relations between the $\sigma^{t,i,a}$'s:
\[ \begin{array}{l}
\sigma^{t+1,j,b} \cdot  \sigma^{t,i,a} =  \sigma^{t+1,i+1,a} \cdot  \sigma^{t,j,b} \\
\hspace{60pt} {\rm for}        \   1 \leq j < i \leq  t  ,   2 \leq {\rm min}( n_{i}, n_{j}), 0 < a < n_{i}, 0 < b < n_{j}  \\
 \\
\sigma^{t+1,i,b} \cdot  \sigma^{t,i,a} =  \sigma^{t+1, i+1, a-b } \cdot \sigma^{t,i,b} \\
 \hspace{60pt} {\rm for}        \   1 \leq j = i \leq  t  ,   2 \leq {\rm min}( n_{i}, a, n_{i} -b ), 0 < a < n_{i}, 0 < b < a  \\
 \\
 \sigma^{t+1,i+1,b} \cdot  \sigma^{t,i,a} =   \sigma^{t+1,i,a }   \cdot  \sigma^{t, i, a + b} \ 0 < b, 2  \leq a,  \ a+b < n_{i} \\
 \\
 \sigma^{t+1,j+1,b} \cdot  \sigma^{t,i,a} = \sigma^{t+1, i , a}  \cdot \sigma^{t,j,b} \  2 \leq {\rm min}( n_{i}, n_{j}), 0 < a < n_{i}, 0 < b < n_{j}, i+2 \leq j.   \\
 \\
 \end{array} \]
 Lastly we have some compositional relations involving the $\tau(K)$'s.
 
We have an obvious notion of one ordered partition of $n$ refining another, for example $\kappa(\alpha)$ refines $\alpha$ which we may denote by $\alpha \prec \kappa(\alpha) $. There is another obvious relation of $\alpha$ being a recombination of $\kappa(\alpha)$ denoted by $\kappa(\alpha) < \alpha$ meaning that $\alpha$ is obtained from $\kappa(\alpha)$ by iterated addition of consecutive integers in the partition.

We have a composition of $\partial^{t,i}$'s 
\[  \partial(\kappa(\alpha) < \alpha) :  \kappa(\alpha)  \longrightarrow \alpha \]
and of $\sigma^{t,i,a}$'s
\[ \sigma(\alpha \prec \kappa(\alpha) ) : \alpha \longrightarrow \kappa(\alpha). \]

The final family of compositional relations is
\[ \begin{array}{l}
  \partial( \kappa(\beta) <  \beta  )  \cdot \tau(K) \cdot  \sigma(\alpha \prec \kappa(\alpha) ) \\
  \\
  =    \sigma(\gamma \prec \beta ) \cdot  \partial( \alpha <  \gamma  ) 
\end{array} \]
where the $\gamma$'s for which the compositions make sense will be illustrated in small examples.

 \section{Small example}
 
 Let us try a small example to see what makes sense. In particular  I want to understand how the $K$ determines the $\gamma$.
 
 In the first small example $K$ will be $2 \times 2$. 
 \[ K = \left( \begin{array}{ccc}
k_{1,1} & k_{1,2}  \\
\\
k_{2,1}&k_{2,2} \\
\end{array} \right)  \]
such that 
\[ \begin{array}{l}
  k_{1,1} + k_{1,2}  = a_{1}   \ {\rm  and}    \   k_{2,1} + k_{2,2}  = a_{2}     \\
  \\
   k_{1,1} + k_{2,1}  = b_{1}   \ {\rm  and}    \   k_{1,2} + k_{2,2}  = b_{2}     \\
  \\
    k_{1,1} + k_{1,2} +  k_{2,1} + k_{2,2}  = a_{1} + a_{2} = b_{1} + b_{2} = n .  \\
    \end{array} \]   
 then the following two compositions are equal
    \[   \begin{array}{c}
     (a_{1}, a_{2})  \\
     \\
 \hspace{80pt}     \downarrow \    \sigma^{3, 3, k_{2,1}} \cdot  \sigma^{2, 1, k_{1,1}} =   \sigma^{3,1, k_{1,1}} \cdot \sigma^{2,2,k_{2,1}} \\
     \\
   ( k_{1,1} , k_{1,2} ,   k_{2,1} , k_{2,2} ) \\
   \\
   \downarrow  \ \tau(K) \\
   \\
    ( k_{1,1} , k_{2,1} ,   k_{1,2} , k_{2,2} )   \\
    \\
    \downarrow \  \partial^{3,1} \cdot  \partial^{4, 3} =   \partial^{3,2} \cdot  \partial^{4, 1}  \\
  \\
    (b_{1}, b_{2}) 
    \end{array}  \]
    and
  \[ \begin{array}{c}
 (a_{1}, a_{2})  \\
 \\
 \downarrow  \  \partial^{2,1}   \\
 \\
 (n) \\
 \\
 \downarrow \ \sigma^{1,1, b_{1}}  \\
 \\
  (b_{1}, b_{2})  
  \end{array}\]
  
  Next suppose that 
  \[ K = \left( \begin{array}{ccc}
k_{1,1} & k_{1,2} &  k_{1,3}  \\
\\
k_{2,1}&k_{2,2}&  k_{2,3} \\
\end{array} \right)  \]
such that 
\[ \begin{array}{l}
k_{1,1} + k_{1,2} +  k_{1,3} = a_{1}, \  k_{2,1} + k_{2,2} +   k_{2,3} = a_{2} \\
\\
k_{1,1} + k_{2,1} = b_{1},  \ k_{1,2} + k_{2,2} = b_{2},  \ k_{1,3} + k_{2,3} = b_{3}   \\
\\
a_{1} + a_{2} = b_{1} + b_{2} + b_{3} = k_{1,1} + k_{1,2} +  k_{1,3} + k_{2,1} + k_{2,2} +   k_{2,3}  = n.
\end{array} \]
The following compositions are equal in this example
\[ \begin{array}{c}
 (a_{1}, a_{2}) \\
 \\
 \downarrow   \sigma^{5,5,k_{2,2}} \sigma^{4,4,k_{2,1}}\sigma^{3,2,k_{1,2}}\sigma^{2,1,k_{1,1}} \\
 \\
 (k_{1,1} , k_{1,2} ,  k_{1,3} ,  k_{2,1} , k_{2,2} ,   k_{2,3}  )   \\
 \\
 \downarrow \ \tau(K)  \\
 \\
 (k_{1,1}, k_{2,1}, k_{1,2}, k_{2,2}, k_{1,3}, k_{2,3}) \\
 \\
 \downarrow     \partial^{4,3}  \partial^{5,2} \partial^{6,1}\\
 \\
 (b_{1}, b_{2}, b_{3}) 
\end{array} \]

\[ \begin{array}{c}
(a_{1}, a_{2}) \\
\\
\downarrow   \   \partial^{2,1} \\
\\
(n) \\
\\
\downarrow \    \sigma^{2, 2, b_{2}}  \sigma^{1,1,b_{1}} \\
\\
(b_{1}, b_{2}, b_{3})
\end{array} \]

If $K$ is the direct sum of two submatrices the relation is a ``product'' of two small compositional relations.
An example follows:
\[ K = \left( \begin{array}{ccccc}
k_{1,1} & k_{1,2} & k_{1,3} & 0 & 0  \\
\\
k_{2,1}&k_{2,2}& k_{2,3}& 0 & 0 \\
\\
0& 0 &  0 &  k_{3,4} &  k_{3,5} \\
\\
0& 0 &  0 &  k_{4,4} &  k_{4,5} \\
\end{array} \right)  \]
such that 
\[ \begin{array}{l}
k_{1,1} + k_{1,2} + k_{1,3} = a_{1}, \ k_{2,1} + k_{2,2} +  k_{2,3} = a_{2} , a_{1} + a_{2} = n_{1}\\
\\
 k_{3,4} +  k_{3,5}=a_{3}, \  k_{4,4} +  k_{4,5}  = a_{4}, \ a_{3} + a_{4} = n_{2}  \\
 \\ 
 k_{1,1} + k_{2,1} = b_{1}, \  k_{1,2} + k_{2,2} = b_{2} , \ k_{1,3} + k_{2,3} = b_{3}, \\
 \\
 k_{3,4} + k_{4,4}  = b_{4}, \ k_{3,5} + k_{4,5} = b_{5} , n = n_{1} + n_{2} \\
\end{array} \]
\[ \begin{array}{c}
(a_{1}, a_{2}, a_{3}, a_{4}) \\
\\
\downarrow  \sigma( (a_{1}, a_{2}) \prec (k_{1,1} \ldots , k_{2,3})) \sigma( (a_{3}, a_{4}) \prec (k_{3,4} \ldots , k_{4,5})) \\
\\
(k_{1,1} , k_{1,2} , k_{1,3}, k_{2,1} , k_{2,2} ,  k_{2,3}, k_{3,4} ,  k_{3,5} , k_{4,4} ,  k_{4,5})
\\
\\
\downarrow \ \tau(K) \\
\\
(k_{1,1} , k_{2,1} , k_{1,2}, k_{2,2} , k_{1,3} ,  k_{2,3}, k_{3,4} ,  k_{4,4} , k_{3,5} ,  k_{4,5})   \\
\\
\downarrow   \  \partial( (k_{1,1}, \ldots , k_{2,3}) <  (b_{1},b_{2}, b_{3})) \partial(  (k_{3,4}, \ldots , k_{4,5}) <  (b_{4},b_{5}) ) \\
\\
(b_{1}, b_{2}, b_{3}, b_{4}, b_{5}) 
\end{array} \]
Observe that $ \partial( (k_{1,1}, \ldots , k_{2,3}) <  (b_{1},b_{2}, b_{3})) \cdot  \partial(  (k_{3,4}, \ldots , k_{4,5}) <  (b_{4},b_{5}) )$ and $ \sigma( (a_{1}, a_{2}) \prec (k_{1,1} \ldots , k_{2,3})) \cdot \sigma( (a_{3}, a_{4}) \prec (k_{3,4} \ldots , k_{4,5})) $ are the ``product'' (or disjoint union) of two disjoint $\partial$'s or $\sigma$'s respectively.

The above composition is equal to the following one:
\[ \begin{array}{c}
(a_{1}, a_{2}, a_{3}, a_{4}) \\
\\
\downarrow  \ \partial((a_{1}, a_{2}) < (n_{1})) \partial((a_{3}, a_{4}) < (n_{2})) \\
\\
(n_{1}, n_{2}) \\
\\
\downarrow   \  \sigma( (n_{1}) \prec (b_{1}, b_{2}, b_{3}))  \sigma( (n_{2}) \prec (b_{4}, b_{5}))\\
\\
(b_{1}, b_{2}, b_{3}, b_{4}, b_{5}) 
\end{array} \]

\section{Realisations}

In the previous section I have given examples of compositional identities in the Hopflike category of \S2.
A Hopflike group is a contravariant from that category to groups.
A familiar example is given \cite{AVZ81} by
\[ A(a_{1}, \ldots , a_{t} ) = R(GL_{a_{1}}{\mathbb F}_{q} \times \ldots \times GL_{a_{t}}{\mathbb F}_{q})  \]
since $A(a_{1}, \ldots , a_{t} )  = \otimes_{i=1}^{t} \ R(GL_{a_{i}}{\mathbb F}_{q} )$ with the $\partial$'s corresponding to iterated multiplication in the Hopf algebras of \cite{AVZ81} and the $\sigma$'s corresponding to the iterated comultiplication, with a ``notional'' copy of $Z$ in degree zero, is the tensor product of countable number of copies of the PSH-algebra given by representations of general linear groups over ${\mathbb F}_{q}$.

I think that 
\[ A_{+}(a_{1}, \ldots , a_{t} ) = R_{+}(GL_{a_{1}}{\mathbb F}_{q} \times \ldots \times GL_{a_{t}}{\mathbb F}_{q})  \]
(notation as in \cite{Sn94} and \cite{Sn88}). This is based upon the partial verification in \cite{Sn20c}.

I also think that 
\[ A_{cmc}(a_{1}, \ldots , a_{t} )  =      {\mathcal H}_{cmc}(GL_{a_{1}}F \times \ldots \times GL_{a_{t}}F)  \]
is another example, where $F$ is a $p$-adic local field and $  {\mathcal H}_{cmc}(G)$ is the hyperHecke algebra of $G$.

\section{Further details concerning $A(a_{1}, \ldots , a_{t})$ }

We have the functor values on objects given by
\[ A(a_{1}, \ldots , a_{t} ) = R(GL_{a_{1}}{\mathbb F}_{q} \times \ldots \times GL_{a_{t}}{\mathbb F}_{q}) . \]
The morphism
 \[   \partial^{t,i} :    ( n_{1}, \ldots , n_{t}) \longrightarrow  ( n_{1}, \ldots , n_{i-1}, n_{i} + n_{i+1}, n_{i+2}, \ldots n_{t}) \] 
 induces a homomorphism of abelian groups
 \[   \begin{array}{c}
  R(GL_{n_{1}}{\mathbb F}_{q} \times \ldots  \times  GL_{n_{i}+n_{i+1}}{\mathbb F}_{q} \times \ldots  GL_{n_{t}}{\mathbb F}_{q} )  \\
 \\
 \downarrow  \    A( \partial^{t,i})  \\
 \\
   R(GL_{n_{1}}{\mathbb F}_{q} \times \ldots  \times  GL_{n_{i}}{\mathbb F}_{q} \times   GL_{n_{i+1}}{\mathbb F}_{q} \times \ldots  GL_{n_{t}}{\mathbb F}_{q} )  
 \end{array}   \]
 sending a representation $V$ to 
 $V^{   1^{ n_{i-1}} \times  U_{n_{i},n_{i+1}}  \times  1^{n_{t} - n_{i+1}} } $. In particular, if 
 \[ V = V_{1} \otimes V_{2} \otimes \ldots \otimes  V_{i-1}  \otimes V' \otimes  V_{i+2} \otimes \ldots \otimes V_{t}  \]
 where $V'$ is a representation of $ GL_{n_{i}+n_{i+1}}{\mathbb F}_{q}$ and $V_{i}$ is a representation of 
 $GL_{n_{i}}{\mathbb F}_{q} $ then 
 \[ A( \partial^{t,i})(V) =  V_{1} \otimes V_{2} \otimes \ldots \otimes  V_{i-1}  \otimes V'^{U_{n_{i}, n_{i+1}} }\otimes  V_{i+2} \otimes \ldots \otimes V_{t} \]
 so 
 \[ A( \partial^{t,i}) = 1^{\otimes n_{i-1}} \otimes  (m^{*})_{n_{i}, n_{i+1}}  \otimes 1^{\otimes n_{t} - n_{i+1}} \]
 which is the external tensor product of a lot of identity maps with the $(n_{i}, n_{i+1})$-component of the PSH-algebra comultiplication.
 
 The map 
  \[ \sigma^{t,i,a} :   ( n_{1}, \ldots , n_{t}) \longrightarrow   ( n_{1}, \ldots , n_{i-1}, a , n_{i} -a,  n_{i+1}, \ldots n_{t}) \]
  induces a homomorphism of abelian groups
  \[ \begin{array}{c}
    R(GL_{n_{1}}{\mathbb F}_{q} \times \ldots  \times  GL_{n_{i-1}}{\mathbb F}_{q} \times  GL_{a}{\mathbb F}_{q} \times   GL_{n_{i} - a}{\mathbb F}_{q} \times \ldots  GL_{n_{t}}{\mathbb F}_{q} )  \\
 \\
 \downarrow  \   A(\sigma^{t,i,a})  \\
 \\
  R(GL_{n_{1}}{\mathbb F}_{q} \times \ldots  \times  GL_{n_{i-1}}{\mathbb F}_{q} \times  GL_{n_{i}}{\mathbb F}_{q} \times \ldots  GL_{n_{t}}{\mathbb F}_{q} )  
  \end{array} \]
which sends a representation $V$ first to its inflation to 
\[ GL_{n_{1}}{\mathbb F}_{q} \times \ldots  \times  GL_{n_{i-1}}{\mathbb F}_{q} \times  P_{a, n_{i} - a} \times \ldots  GL_{n_{t}}{\mathbb F}_{q} \]
and then induces it up to 
\[ GL_{n_{1}}{\mathbb F}_{q} \times \ldots  \times  GL_{n_{i-1}}{\mathbb F}_{q} \times  GL_{n_{i}}{\mathbb F}_{q} \times \ldots  GL_{n_{t}}{\mathbb F}_{q}  . \]
If 
\[ V = V_{1} \otimes V_{2} \otimes \ldots \otimes  V_{i-1}  \otimes V_{a}'  \otimes  V_{n_{i} -a}' \otimes  V_{i+2} \otimes \ldots \otimes V_{t}  \]
it is sent to 
\[ V_{1} \otimes V_{2} \otimes \ldots \otimes  V_{i-1}  \otimes m(V_{a}'  \otimes  V_{n_{i} -a}') \otimes  V_{i+2} \otimes \ldots \otimes V_{t} \]
where $m$ is the multiplication in the PSH-algebra.

An ordered partition has two significant integers relating to it. These are the partition sum of 
$(a_{1}, \ldots , a_{t})$ given by $a_{1} + \ldots + a_{t} = n$ and the partition length, which equals $t$.
The basic morphisms make the following changes:
\[ \begin{array}{|c|c|c|}
\hline 
\hline
{\rm morphism} & {\rm partition \ sum} & {\rm partition \ length} \\
\hline
\hline
\partial^{t,i} & 0 & -1   \\
\hline
\sigma^{t,i,a} & 0 & + 1 \\
\hline
\tau(K) & 0 & 0 \\
\hline
\hline
\end{array} \]

The PSH-algebra of (\cite{AVZ81} Chapter III) is given by $R = \oplus_{n \geq 0} \ R(GL_{n}{\mathbb F}_{q})$
which corresponds to the products of $GL_{n}$'s on length $1$ (with a fake convention for $GL_{0}$).The product has components of the form 
\[ R(GL_{n}{\mathbb F}_{q}) \otimes R(GL_{m}{\mathbb F}_{q}) = R(GL_{n}{\mathbb F}_{q}) \times GL_{m}{\mathbb F}_{q}) \longrightarrow  R(GL_{n+m}{\mathbb F}_{q}) \]
and the coproduct goes to the sum of length two partitions with the same sum
\[ R(GL_{n}{\mathbb F}_{q}) \longrightarrow  \oplus_{{\rm length \ 2, \ sum \  n}} \ R(GL_{a_{1}}{\mathbb F}_{q}
\times GL_{a_{2}}{\mathbb F}_{q} ) .  \]

The full $(n_{1}, n_{2})$-component of the coproduct may be described as
\[  A(\partial^{2,1}) :  R(GL_{n}{\mathbb F}_{q})  \longrightarrow   R(GL_{n_{1}}{\mathbb F}_{q} \times  
GL_{n_{2}}{\mathbb F}_{q} )   \]
providing that we interpret $ A(\partial^{2,0}) $ and $ A(\partial^{2,2})$ as the identity map. 

The product
\[ \oplus_{{\rm length \ 2, \ sum \  n}} \ R(GL_{a_{1}}{\mathbb F}_{q}
\times GL_{a_{2}}{\mathbb F}_{q} )  \longrightarrow   R(GL_{n}{\mathbb F}_{q})  \]
has $(n_{1}, n_{2})$-component given by
\[ A(\sigma^{1, 1, n_{1}}) :   R(GL_{n_{1}}{\mathbb F}_{q}
\times GL_{n_{2}}{\mathbb F}_{q} )  \longrightarrow   R(GL_{n}{\mathbb F}_{q})  \]
providing that we interpret $ A(\sigma^{t,0,0}) $ and $ A(\sigma^{t,t,0})$ as the identity map. 

The $r \times s$ matrix of (strictly) positive integers
\[ K = \left( \begin{array}{ccccc}
k_{1,1} & k_{1,2} & \ldots & \ldots & k_{1,s}  \\
\\
k_{2,1}&k_{2,2}& \ldots & \ldots & k_{2,s} \\
\\
\vdots & \vdots &  \vdots &  \vdots &  \vdots \\
\\
k_{r,1} & k_{r,2} & \ldots & \ldots & k_{r,s}, \\
\end{array} \right)  \]
such that $ \sum_{j} \ k_{i,j} = a_{i}$ for $1 \leq i \leq r$ and  $ \sum_{i} \ k_{i,j} = b_{j}$ for $1 \leq j \leq s$, , from the viewpoint of ordered products $GL_{n}$'s, sends $GL_{n}$'s strung out by rows
\[ GL_{k_{1,1}}{\mathbb F}_{q} \times GL_{k_{1,2}}{\mathbb F}_{q} \times GL_{k_{1,3}}{\mathbb F}_{q}
\ldots   \ldots   \ldots   GL_{k_{r,s,1}}{\mathbb F}_{q}  \]
to the same groups but in an order product strung out by columns
\[ GL_{k_{1,1}}{\mathbb F}_{q} \times GL_{k_{2,}}{\mathbb F}_{q} \times GL_{k_{3,1}}{\mathbb F}_{q}
\ldots  \ldots   \ldots   GL_{k_{r,s,1}}{\mathbb F}_{q} . \]
Hence $\tau(K)$ preserves both partition sum and partition length.

In this example the verification of the PSH-algebra axioms in \cite{AVZ81} shows that any two compositions of the natural transformations obtained by applying $A(-)$ to the generating morphism in the Hopflike category 
with the same range and domain, indexed by ordered partitions, are equal. Furthermore, with the Schur inner product, this example is self-dual (in the sense that $m^{*}$ is the adjoint of $m$).

Next I want to do the first of the other "products" in this example. This has the form
\[ m: A(a_{1}) \otimes A(b_{1}, \ldots , b_{t}) \longrightarrow  \oplus_{i=1}^{t} \ A(b_{1}, \ldots ,b_{i-1},  a_{1} + b_{i},
b_{i+1}, \ldots , b_{t})  \]
given by
\[ m(V \otimes (V_{1} \otimes \ldots  \otimes V_{t})) = \oplus_{i=1}^{t} \ V_{1} \otimes V_{i-1} \otimes 
 {\rm IndInf}_{G_{a_{1} } \times G_{b_{i}}}^{G_{a_{1} + b_{i}}}(V \otimes V_{i}) \otimes \ldots V_{i+1} \otimes \ldots \otimes V_{t} .\]
 The neater way to think about this, which generalises to all pairs of ordered partitions is to add zeroes in all possible ways to make the partitions the same length, namely
 \[ (a_{1}, 0, \ldots , 0), (0, a_{1}, 0 , \ldots ,0), \ldots , (0,0, \ldots , 0, a_{1}) \]
 and then add up each of the "products" of the form
 \[   \begin{array}{l}
  (W_{1} \otimes \ldots  \otimes W_{t})  \cdot  (V_{1} \otimes \ldots  \otimes V_{t}) ) \\
 \\
 =   {\rm IndInf}_{G_{a_{1} } \times G_{b_{i}}}^{G_{a_{1} + b_{1}}}(W_{1} \otimes V_{1}) \otimes \ldots 
 \otimes  {\rm IndInf}_{G_{a_{t} } \times G_{b_{t}}}^{G_{a_{t} + b_{t}}}(W_{t} \otimes V_{t}) 
 \end{array}  \]
 which makes sense because, for example, if $a_{i}=0$ then  $G_{a_{i}} = \{ 1 \}$ and $W_{i} = 1$ so that
 $ {\rm IndInf}_{G_{a_{i} } \times G_{b_{i}}}^{G_{a_{i} + b_{i}}}(W_{i} \otimes V_{i}) = V_{i}$ .
 
 This procedure gives an associative (because the inflation-induction is associative) product on 
 \[   A = \oplus_{(a_{1}, \ldots , a_{t} )}  \ A(a_{1}, \ldots , a_{t} ) \]
 where the sum is over the equivalence classes of ordered partitions under the ``adding zeroes'' equivalence relation (see the earlier footnote).
 
 The coproduct on $A$ must surely be
 \[ A(n_{1}, \ldots , n_{t}) \longrightarrow  \oplus_{i=1}^{t} \oplus_{a=0}^{n_{i}} \ A(n_{1}, \ldots ,n_{i-1}, a, n_{i}-a, \ldots , n_{t}) \]
 sending $V_{1} \otimes  \ldots \otimes V_{t}$ to the direct sum whose $(i,a)$-component is 
 \[V_{1} \otimes \ldots \otimes V_{i-1} \otimes m^{*}_{a, n_{i}-a}(V_{1}) \otimes \ldots  \otimes V_{t} . \]
 
Recall that in a Hopf algebra $H$ we have a commutative diagram which states that the following two homomorphisms are equal\footnote{In the following discussion one has to try to remember to distinguish between $m^{*}$ the coproduct in the PSH-algebra of \cite{AVZ81} and the coproduct in $A$ which I am trying to define!}
\[ H \otimes  H  \stackrel{m}{\longrightarrow}  H   \stackrel{m^{*}}{\longrightarrow}  H \otimes  H \]
and
\[  H \otimes  H  \stackrel{m^{*} \otimes m^{*}}{\longrightarrow} H \otimes  H  \otimes H \otimes  H   
\stackrel{1 \otimes T \otimes 1}{\longrightarrow }   H \otimes  H  \otimes H \otimes  H    
\stackrel{m \otimes m}{\longrightarrow}  H \otimes H . \]

We know that the analogue of this condition holds if we start with $A(a_{1}) \otimes A(b_{1})$. Let us examine what happens if we start with $A(a_{1}) \otimes A(b_{1}, b_{2})$ and without any loss of generality we may start with the images of $W  \otimes V_{1} \otimes V_{2}$. Under the product this goes to the pair
$(W \otimes 1 \otimes V_{1} \otimes V_{2}, 1 \otimes W \otimes V_{1} \otimes V_{2})$ and then to 
\[   (V_{1} \otimes {\rm IndInf}(W \otimes V_{2}) ,   {\rm IndInf}(W \otimes V_{1}) \otimes V_{2}) \in 
A(b_{1}, a_{1} + b_{2}) \oplus  A(a_{1} + b_{1}, b_{2})  .\]
Next it goes to 
\[ \begin{array}{l}
m^{*}(V_{1}) \otimes  {\rm IndInf}(W \otimes V_{2})  \  \oplus \   V_{1} \otimes 
m^{*}( {\rm IndInf}(W \otimes V_{2})) \\
\\
\oplus  m^{*}( {\rm IndInf}(W \otimes V_{1})) \otimes V_{2}  \   \oplus  
{\rm IndInf}(W \otimes V_{1}) \otimes  m^{*}(V_{2}) 
\end{array} \]

Going via the lower route $W  \otimes V_{1} \otimes V_{2}$ goes first to
\[  m^{*}(W)  \otimes m^{*}(V_{1} \otimes V_{2}) = m^{*}(W) \otimes m^{*}(V_{1}) \otimes V_{2}  \oplus 
m^{*}(W) \otimes V_{1} \otimes m^{*}(V_{2}). \]
 The tensor-factors of this are indexed by
partitions  $(x, a_{1} - x, y, b_{1}-y, b_{2})$ and $(x, a_{1} - x, b_{1}, y' , b_{2}-y')$. If $U_{x,a_{1}-x} \subset GL_{a_{1}}$,  $U_{y, b_{1}-y} \subset GL_{b_{1}}$, $U_{y', b_{2} - y'} \subset GL_{b_{2}}$ are usual the unipotent subgroups.

Let us consider the summand  $m^{*}(W) \otimes m^{*}(V_{1}) \otimes V_{2}$. $W \otimes V_{1} \otimes V_{2}$ is a representation of $GL_{a_{1}} \times GL_{b_{1}} \times GL_{b_{2}}$ and applying $m^{*} \otimes m^{*}$ gives
  \[    W^{U_{x,a_{1}-x}} \otimes  (V_{1})^{U_{y,b_{1}-x}} \otimes V_{2}   \]
 as a representation of $GL_{x} \times GL_{a_{1}-x} \times GL_{y} \times GL_{b_{1}-y} \times GL_{b_{2}}$
 then $1 \times T \times 1$ switches the product of $GL$'s to 
 $GL_{x}  \times GL_{y}  \times GL_{a_{1}-x}  \times GL_{b_{1}-y} \times GL_{b_{2}}$. To apply the Hopflike $m \otimes m$ to this element of $A(x,y) \otimes A(a_{1}-x, b_{1}-y, b_{2})$ we muust add in the zero in three places to give
 \[ (x,y,0; a_{1}-x, b_{1}-y, b_{2}) , \  (x,0, y); a_{1}-x, b_{1}-y, b_{2}) \  (0, x,y); a_{1}-x, b_{1}-y, b_{2}) \]
 then we apply the Hopf algebra $m \otimes m \otimes m$. The first of the three $6$-tuples sums to give
 $m^{*}m(W \otimes V_{1}) \otimes V_{2}$ which is the bottom left term of the quartet of terms for the upper route.

\section{Further details concerning $A_{+}(a_{1}, \ldots , a_{t})$}

 \section{Unexplained Appendix from \cite{Sn20c}: The general $m$ and $m^{*}$ condition}

Suppose that we have two ordered partitions of $n$ (c.f. \cite{AVZ81} p. 170) 
\[  \alpha = (a_{1}, a_{2}, \ldots , a_{r} ) \ {\rm and}  \  \beta = (b_{1}, b_{2}, \ldots , b_{s})  . \]
Suppose that we have an $r \times s$ matrix of non-negative integers 
\[ K = \left( \begin{array}{ccccc}
k_{1,1} & k_{1,2} & \ldots & \ldots & k_{1,s}  \\
\\
k_{2,1}&k_{2,2}& \ldots & \ldots & k_{2,s} \\
\\
\vdots & \vdots &  \vdots &  \vdots &  \vdots \\
\\
k_{r,1} & k_{r,2} & \ldots & \ldots & k_{r,s} \\
\end{array} \right)  \]
such that $ \sum_{j} \ k_{i,j} = a_{i}$ for $1 \leq i \leq r$ and  $ \sum_{i} \ k_{i,j} = b_{j}$ for $1 \leq j \leq s$.
Therefore we have ordered partitions of $n$
\[ \kappa(\alpha) :  (k_{1,1}, k_{1,2}, \ldots , k_{1,s}, k_{2,1}, \ldots , k_{2,s}, \ldots , k_{r,1}, \ldots , k_{r,s}  ) \]
and
\[ \kappa(\beta) :( k_{1,1}, k_{2,1}, \ldots , k_{r,1}, k_{1,2}, \ldots , k_{r,2}, \ldots , k_{1,s}, \ldots , k_{r,s}   ) . \]
Now think of a partition as subdividing the ordered set $(1,2,3, \ldots , n)$ into abutting intervals (called ``blocks'' in \cite{AVZ81}) of which, for $\alpha$ say, the first is $I_{1} = (1,2, \ldots , a_{1})$, the second is $I_{2} = (a_{1}+1, a_{1}+2, \ldots , a_{1}+a_{2})$ and so on up to $I_{r}$. Similarly the intervals for 
$\beta$ are $J_{1}, \ldots , J_{s}$. Then in the two partitions, $\kappa(\alpha)$ and $ \kappa(\beta)$, involving the $k_{i,j}$'s given above
there is a permutation $\sigma_{K} \in \Sigma_{n}$ which shifts (by a translation) the $k_{i,j}$-interval
of the $\kappa(\alpha)$ into the $i$-th interval of $J_{j}$ in $\kappa(\beta)$. For example, the integers in the $k_{1,2}$-interval in the upper partition goes, by a linear shift, to the interval of integers $(b_{1}+1, b_{1}+2, \ldots , b_{1}+k_{1,2})$ in the lower partition.

We have an obvious notion of one ordered partition of $n$ refining another, for example $\kappa(\alpha)$ refines $\alpha$ which we may denote by $\alpha \prec \kappa(\alpha) $. There is another obvious relation of $\alpha$ being a recombination of $\kappa(\alpha)$ denoted by $\kappa(\alpha) < \alpha$ meaning that $\alpha$ is obtained from $\kappa(\alpha)$ by iterated addition of consecutive integers in the partition.

As in \cite{Sn20b} let ${\mathcal H}_{cmc}(G)$ denote the hyperHecke algebra of a locally $p$-adic Lie group $G$. For an ordered partition $\alpha$ write 
$GL(\alpha) = GL_{a_{1}}K \times GL_{a_{2}}K \times  \ldots \times  GL_{a_{r}}K $ and consider  ${\mathcal H}_{cmc}(GL(\alpha))$.

I believe that for each refinement $\alpha \prec \kappa(\alpha) $ we have a ``coproduct''
\[ m^{*}(\alpha \prec \kappa(\alpha)) : {\mathcal H}_{cmc}(GL(\alpha)) \longrightarrow 
{\mathcal H}_{cmc}(GL(\kappa(\alpha))) \]
and for each recombination we have a ``multiplication''
\[ m(\kappa(\alpha) < \alpha) : {\mathcal H}_{cmc}(GL(\kappa(\alpha))) \longrightarrow  {\mathcal H}_{cmc}(GL(\alpha)).  \]
 I expect each of these maps to be associative. Furthermore for each $\alpha$ and $\beta$ are above 
 I expect the following compositions to be equal:
 \[ \begin{array}{c}
   {\mathcal H}_{cmc}(GL(\alpha))    \\
  \\
\hspace{60pt}  \downarrow \  m^{*}(\alpha \prec \kappa(\alpha)) \\
  \\
  {\mathcal H}_{cmc}(GL(\kappa(\alpha))) \\
  \\
\hspace{60pt}     \downarrow \  (\sigma_{K})_{*} \\
    \\
      {\mathcal H}_{cmc}(GL(\kappa(\beta))) \\
      \\
  \hspace{60pt}     \downarrow \ m(\kappa(\beta) < \beta)   \\
      \\
        {\mathcal H}_{cmc}(GL(\beta))
  \end{array} \]
 
 \[ \begin{array}{l}
    {\mathcal H}_{cmc}(GL(\alpha))    \\
  \\
\hspace{60pt}  \downarrow \  m(\alpha <  \gamma) \\
\\
  {\mathcal H}_{cmc}(GL(\gamma))    \\
  \\
  \hspace{60pt}  \downarrow \  m^{*}(\gamma \prec \beta) \\
\\
  {\mathcal H}_{cmc}(GL(\beta))    \\
 \end{array} \]
 for any $\gamma$ for which these compositions make sense.
 
 My evidence for these beliefs comes from the analogue of Hopf algebra identity for the hyperHecke alagebras when the ordered partitions $\alpha$ and $\beta$ have length two.
 
 \section{ Low dimensional combinatorics of ``Hopflike'' algebras}
 
In this section an ordered partition $\alpha = (a_{1}, \ldots , a_{r})$ will allow the $a_{i}$'s to be non-negative integers, $a_{i} \geq 0$. Suppose we have abelian groups $R(a_{i})$ which are going to imitate the case of \cite{AVZ81} in which $R(a_{i})$ equals the complex representation ring of $GL_{a_{i}}{\mathbb F}_{q}$. In this case we have $R(0) = {\mathbb Z}$ and we have Hopf algebra structures
\[ m: R(a) \otimes R(b) \longrightarrow R(a+b) \]
and 
\[ m^{*} : R(a) \longrightarrow  \oplus_{u=0}^{a} \ R(u) \otimes R(a-u) . \] 
The Hopf algebra condition asserts that the composite homomorphism
\[ R(a) \otimes R(b) \stackrel{m}{\longrightarrow}  R(a+b)  \stackrel{m^{*}}{\longrightarrow}  
  \oplus_{j=0}^{a+b}   \ R(j) \otimes R(a+b-j)  \]
   is equal to the composite homomorphism 
\[ \begin{array}{l}
R(a) \otimes R(b)  \stackrel{{\rm shuffle} (m^{*} \otimes m^{*})}{\longrightarrow} 
\oplus_{u=0}^{a} \oplus_{v=0}^{b} R(u) \otimes R(v) \otimes R(a-u) \otimes R(b-v) \\
\\
\hspace{50pt}  \stackrel{(m \otimes m)}{\longrightarrow}  \oplus_{u,v}  \ R(u+v) \otimes R(a+b-u-v).
\end{array}  \]

The combinatorics of the $u,v, j$'s is summed up in the following rectangle.
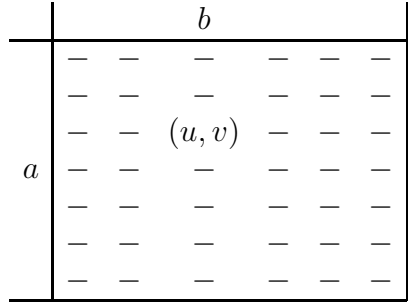
\[ \begin{array}{c|cccccc|}
&&&b&&& \\
\hline
&-& - & - & - & - & - \\
&-& - & - & - & - & - \\
&-& - & (u,v) & - & - & - \\
a&-& - & - & - & - & - \\
&-& - & - & - & - & - \\
&-& - & - & - & - & - \\
&-& - & - & - & - & - \\
\hline
\end{array}  \]

Suppose we have $R(a) \otimes (R(b) \otimes R(c))$\footnote{See Question \ref{8.1}: I have also to do $((R(a) \otimes R(b)) \otimes R(c)$ and $(R(a)\otimes - \otimes R(c)) \otimes R(b)$ in full combinatorial detail. In particular, equivariance with respect to permutations.} we define a multiplication $m$ on this triple tensor product, which has bidegree $(1,2)$ in terms of lengths of partitions.
\[    \begin{array}{l}
m= 0 :  R(a) \otimes (R(b) \otimes R(c)) \longrightarrow  0 , \ {\rm if} \ a > 0, b > 0, c > 0  \\ 
\\
m = m :  R(a) \otimes (R(b) \otimes R(c)) \longrightarrow   R(b+c) , \ {\rm if} \ a =0, \\
\\ 
m = m :  R(a) \otimes (R(b) \otimes R(c)) \longrightarrow   R(a+c) , \ {\rm if} \ b =0, \\
\\
m = m :  R(a) \otimes (R(b) \otimes R(c)) \longrightarrow   R(a+b) , \ {\rm if} \ c =0. \\
\end{array} \]

I am going to show you the very neat combinatorial formula for the homomorphism
\[ (m \otimes m)  \cdot ({\rm shuffle})(m^{*} \otimes m^{*} \otimes m^{*}),  \]
the triple tensor of the comultiplication following by the double tensor of my $m$, minus
the homoorphism
\[ m^{*} \cdot m , \]
the composition of the Hopf algebra comultiplication following the ``multiplication'' $m$ which I just defined.

When $abc=0$ the two homomorphisms are equal by the Hopf algebra axiom.

When $a > 0, b > 0, c > 0$ the second of these homomorphisms is zero. The $m \otimes m$ in the second map
\[ \oplus_{u,v,w} \ (R(u) \otimes R(v) \otimes R(w)) \otimes (R(a-u) \otimes R(b-v) \otimes R(c-w)) \longrightarrow  ?? \]
is non-zero in the six cases
\[ \begin{array}{l}
u = 0 , b = v \\
\\
u = 0, c = w \\
\\
v=0,  a = u \\
\\
v=0, c = w \\
\\
w=0, a = u \\
\\
w=0,  b = v. \\
\end{array} \]
Assuming commutativity of the Hopf algebra multiplication we obtain for the difference in homomorphisms
\[ \begin{array}{l}
(1) \hspace{30pt} \sum_{w=0}^{c} \ m(R(a) \otimes R(w)) \otimes m(R(b) \otimes R(c-w))  \\
\\
(2) \hspace{30pt}  + \sum_{w=0}^{c} \ m(R(b) \otimes R(w)) \otimes m(R(a) \otimes R(c-w))   \\
\\
(3) \hspace{30pt}  \sum_{v=0}^{b} \ m(R(a) \otimes R(v)) \otimes m(R(c) \otimes R(b-v))  \\
\\
(4) \hspace{30pt}  + \sum_{v=0}^{b} \ m(R(c) \otimes R(v)) \otimes m(R(a) \otimes R(b-v))   \\
\\
(5) \hspace{30pt}  \sum_{u=0}^{a} \ m(R(b) \otimes R(u)) \otimes m(R(c) \otimes R(a-u))  \\
\\
(6) \hspace{30pt}  + \sum_{u=0}^{a} \ m(R(c) \otimes R(u)) \otimes m(R(b) \otimes R(a-u))  .
\end{array} \]

Here comes, without assuming commutativity of the Hopf algebra multiplication, the first stage in the formula for extra terms occasioned by 
\newpage
$R(a) \otimes (R(b) \otimes R(c))$.
\[ \begin{array}{l}
(1) \hspace{30pt} \sum_{w=0}^{c} \ R(a)  \otimes R(0) \otimes R(w)  \otimes R(0) \otimes  R(b) \otimes R(c-w)  \\
\\
(2) \hspace{30pt}  + \sum_{w=0}^{c} \ R(0) \otimes R(b) \otimes R(w) \otimes R(a)  \otimes R(0) \otimes R(c-w)   \\
\\
(3) \hspace{30pt}  \sum_{v=0}^{b} \ R(a) \otimes R(v) \otimes  R(0) \otimes R(0)  \otimes R(b-v)  \otimes R(c) \\
\\
(4) \hspace{30pt}  + \sum_{v=0}^{b} \ R(0) \otimes R(v)  \otimes R(c)  \otimes R(a) \otimes R(b-v)  \otimes R(0) \\
\\
(5) \hspace{30pt}  \sum_{u=0}^{a} \  R(u)  \otimes R(b) \otimes R(0)  \otimes  R(a-u)  \otimes R(0) \otimes R(c)   \\
\\
(6) \hspace{30pt}  + \sum_{u=0}^{a} \ R(u) \otimes R(0) \otimes R(c) \otimes R(a-u) \otimes R(b) \otimes R(0)  .
\end{array} \]

The formula which I am assembling slowly is the sum of six maps
\[ R(a) \otimes (R(b) \otimes R(c))  \mapsto   (n) \ ({\rm for } \ n=1,2,3,4,5,6. \]
Here are the six formulae, each of which is applied to $R(a) \otimes (R(b) \otimes R(c))$\footnote{The permutation given by the middle bracket should be written in the notation of \cite{AVZ81}}.
\[ \begin{array}{l}
(1) \hspace{30pt} \sum_{w=0}^{c} \ (m_{a,w} \otimes m_{b, c-w}) 
( 1 \otimes  \left(  \begin{array}{cc}
0 & I_{w} \\
\\
I_{b} & 0  \\
\end{array} \right) \otimes 1 ) ( 1 \otimes \ 1 \otimes m^{*}_{c})   \\
\\
(2) \hspace{30pt}  +   \sum_{w=0}^{c} \  ( m_{b,w} \otimes m_{a,c-w}) 
( 1 \otimes  \left(  \begin{array}{cc}
0 & I_{w} \\
\\
I_{a} & 0  \\
\end{array} \right) \otimes 1 ) 
(  \left(  \begin{array}{cc}
0 & I_{a} \\
\\
I_{b} & 0  \\
\end{array} \right) \otimes m^{*}_{c} )  \\
\\
(3) \hspace{30pt}  +   \sum_{v=0}^{b} \ ( m_{a,v} \otimes m_{c,b-v}) 
( 1 \otimes  \left(  \begin{array}{cc}
0 & I_{c} \\
\\
I_{b-v} & 0  \\
\end{array} \right) \otimes 1 ) (1 \otimes m^{*}_{b} \otimes 1)  \\
\\
(4) \hspace{30pt}  +   \sum_{v=0}^{b} \  ( m_{a,b-v} \otimes m_{v,a}) 
( 1 \otimes  \left(  \begin{array}{cc}
0 & I_{v} \\
\\
I_{b-v} & 0  \\
\end{array} \right) \otimes 1 ) \\
\\
\hspace{140pt}  ( 1 \otimes m^{*}_{b} \otimes 1)  \left( \begin{array}{ccc} 
0 & 0 & I_{a} \\
\\
0 & I_{b} & 0 \\
\\
I_{c}& 0 & 0 
\end{array}  \right)   \\
\\
(5) \hspace{30pt}  +   \sum_{u=0}^{a} \  ( m_{u,b} \otimes m_{a-u,c})
( 1 \otimes  \left(  \begin{array}{cc}
0 & I_{b} \\
\\
I_{a-u} & 0  \\
\end{array} \right) \otimes 1 ) ( m^{*}_{a} \otimes 1 \otimes 1)  \\
\\
(6) \hspace{30pt}  +  \sum_{u=0}^{a} \  ( m_{u,c} \otimes m_{a-u,b})
(1 \otimes  \left( \begin{array}{ccc}
0 & 0 & I_{c} \\
\\
I_{a-u} & 0 & 0 \\
\\
0 & I_{b} & 0 \\
\end{array} \right) )  ( m^{*}_{a} \otimes 1 \otimes 1)  \\
\\
\end{array} \]

\begin{question}{$_{}$}
\label{8.1}
\begin{em}

Is this expression familiar from $A_{\infty}$-algebras or something similar?

What about adding this term to the other homomorphism to get an identity\footnote{Finding how to modify $m$ in bidegrees $(1,2)$ and $(2,1)$ is reminiscent of my introduction to homotopy coherence - H-spaces $X$ which, say, are homotopy associative by homotopies which fit into a higher dimensional picture of homotopies and so on and on. These ideas developed into operads, and $A_{\infty}$, $E_{\infty}$ spaces, spectra and rings. Around 1970 papers by Jim Stasheff showed how much of the classifying space of $X$ one could make from how much coherence. For example, he started with the projective plane of $X$. The extension of $m$ to bidegree $(1,2)$ and $(2,1)$ is perhaps analogous to getting to the projective plane - a first step.} ?
\end{em}
\end{question}

{\bf Note to myself:}  Once I get my mind straight again, I should forget the algebraic simplification ``Assuming commutativity of the Hopf algebra multiplication'' and write out the extra terms with the correct permutations included everywhere. Then I should add these terms into a formula for a new comultiplication which fixes the Hopf condition. My guess is that these augmented comultiplications generalise to all bidegrees from the cases $(1,2)$ and $(2,1)$. If this sort of Hopf algebra works it should dove-tail nicely with the Langlands `unique irreducible quotient" \cite{AJS78} result and maybe even the modular case of Florian Herzig and et al \cite{AHHV2017}.


\begin{thebibliography}{AVZ81}
 
 \bibitem{AHHV2017}  N. Abe, G. Henniart, F. Herzig and M-F. Vigneras: A classification of irreducible admissible mod $p$ representations of $p$-adic reductive groups; J.A.M.Soc. 30 \#2 (2017) 495-559.
 (http://dx.doi.org/ 10.1090/jams/862.)
 
 \bibitem{AKdS13}  Eran Assaf, David Kazhdan and Ehud de Shalit: Kirillov models and the Breuil-Schneider conjecture for $GL_{2}(F)$; arXiv:1302.3060.2013.

\bibitem{BZ76}  J. Bernstein and A. Zelevinski: Representations of the group $GL(n,F)$ where $F$ is a local non-Archimedean field; Uspekhi Mat. Nauk. {\bf 31} 3 (1976) 5-70.

\bibitem{BZ77}  J. Bernstein and A. Zelevinski: Induced representations of reductive $p$-adic groups I; Ann. ENS {\bf 10} (1977) 441-472.

\bibitem{Bol01} R. Boltje: Monomial resolutions; J.Alg. 246 (2001) 811-848.

\bibitem{Bour58}  N. Bourbaki:  {\em Alg\`{e}bre}; Hermann  (1958).

\bibitem{Bour63}   N. Bourbaki:  {\em  M\'{e}sures de Haar};  Hermann  (1963).

\bibitem{Bour67}   N. Bourbaki:  {\em Vari\'{e}t\'{e}s  diff\'{e}rentielles et analytiques};  Hermann  (1967).

\bibitem{Bour68}   N. Bourbaki:  {\em  Groupes et alg\`{e}bres de Lie};  Hermann  (1968). 

 \bibitem{NB68} N. Bourbaki: Groupes er alg\`{e}bres de Lie; Chapters IV-VI  (1968) Hermann Paris.



\bibitem{FBr56} F. Bruhat: Sur les repr\'{e}sentations induites des groupes de Lie; Bull. Soc. Math. France 84 (1956) 97-205.

\bibitem{FBr61}  F. Bruhat: Distributions sur un groupe localement compact et applications \`{a} l'\'{e}tude des repr\'{e}sentations des groupes $p$-adiques; ; Bull. Soc. Math. France 89 (1961) 43-75.

\bibitem{BH06}  C.J. Bushnell and G. Henniart: {\em The Local Langlands Conjecture for $GL(2)$}; Grund. Math. Wiss. \#335; Springer Verlag (2006).

  
\bibitem{DB96} Daniel Bump: {\em Automorphic forms and representations}; Cambridge studies in advanced math. {\bf 55} (1998).

\bibitem{CEGGPS14} Ana Cariani, Matthew Emerton, Toby Gee, David Geraghty, Vyautas Paskunas and Sug Woo Shin: Patching and the $p$-adic Langlands correspondence; Cambridge J. Math (2014).

\bibitem{HDHH35}  H. Davenport and H. Hasse: Die Nullstellen der Kongruenzzetafunktionen in gewissen zyklischen F�llen; Journal f�r die Reine und Angewandte Mathematik {\bf 172} (1935) 151�182.

  
\bibitem{PD84}  P. Deligne: Le ``centre'' de Bernstein;  {\em Repr\'{e}sentations des groupes r\'{e}ductifs sur un corps local}  Travaux en cours, Hermann, Paris (1984) 1-32.
 
 \bibitem{JAG55} J.A. Green: The characters of the finite general linear groups; Trans. Amer. Math. Soc. 80 (1955) 402-447.
 
 \bibitem{FH2011 }  Florian Herzig: The classification of irreducible admissible mod $p$ representations of a $p$-adic $GL_{n}$; Inv. Math.  {\bf 186} (2011) 373-434.
 
 
 
 \bibitem{HS71}  P.J. Hilton and U. Stammbach: {\em A Course in Homological Algebra}; GTM \#4  (1971)  Springer Verlag.

 
 \bibitem{GDJ76} G.D. James: The irreducible representations of the symmetric groups; Bull. London Math. Soc.  8 (1976) 229-232.
 
 \bibitem{GDJ} G. D. James: {\rm The Representation Theory of the Symmetric Groups}; Springer Verlag Lecture Notes in Math. \#682. 
 
 \bibitem{KMS2011} J. P. Keating, F. Mezzadri and B. Singphu: Rate of convergence of linear functions on the unitary group; J. Phys. A. Math. Theor. 44 (2011) o35204.
 
 \bibitem{Kon63}  T. Kondo: On Gaussian sums attached to the general linear groups over finite fields; J. Math. Soc. Japan vol.15  \#3 (1963) 244-255.

  
\bibitem{IGM80}  I.G. Macdonald: Zeta functions attached to finite general linear groups; Math. Annalen 249 (1980) 1-15.

\bibitem{MSZ89}  J.P. May, V.P. Snaith  and P. Zelewski: A further generalisation of the 
Segal conjecture; Quart. J. Math. Oxford (2) 40  (1989) 457-473.  

  \bibitem{MM65} J.W. Milnor and J.C. Moore: On the structure of Hopf algebras; Annals of Math. (2) 81 (1965) 211-264.

\bibitem{MZ55} D. Montgomery and L. Zippin: {\em Topological Transformation Groups}; Interscience New York (1955).

\bibitem{JMOS18} J.M. O'Sullivan and C.N. Harrison: Myelofibrosis: Clinicopathologic Features, Prognosis and Management; Clinical Advances in Haematology and Oncology {\bf 16} (2) February 2018.1.

\bibitem{AMR2000}  A. M. Robert: {\em A Course in $p$-Adic Analysis}; Grad Texts in Math. \#198, Springer Verlage (2000). 

\bibitem{Ser77}   J-P. Serre: {\em  Linear Representations of Finite Groups }; Grad. Texts in Math.  \# 42  (1977)  Springer-Verlag.

\bibitem{JPS2004} J-P. Serre: Compl\`{e}te R\'{e}ductibilit\'{e}; S\'{e}minaire BOURBAKI 56i\`{e}me
ann\'{e}e 2003-2004 \#932  p. 195 \`{a} 217.

\bibitem{Sh76}  T. Shintani: Two remarks on irreducible characters of finite general linear groups;  J. Math. Soc. Japan {\bf 28} (1976) 396-414.

\bibitem{AJS78} Allan J. Silberger: The Langlands quotient theorem for p-adic groups; Math. Annalen 236 no. 2 (1978) 95-104.

\bibitem{AJS79}  Allan J. Silberger: {\em Introduction to harmonic analysis on $p$-adic reductive groups}; Math. Notes Princeton Unviersity Press (1979).

\bibitem{Sn88} V.P. Snaith: Explicit Brauer Induction; Inventiones Math. 94 (1988) 455-478.

\bibitem{Sn89}  V.P. Snaith: {\em  Topological Methods in Galois Representation Theory}, C.M.Soc Monographs, Wiley (1989) (republished by Dover in 2013).


\bibitem{Sn94}  V.P. Snaith:  {\em Explicit Brauer Induction (with applications to algebra and number theory)}; Cambridge studies in advanced mathematics \#40,Cambridge University Press (1994).
  
\bibitem{Sn18} V.P. Snaith:  {\em Derived Langlands}; World Scientific (2018).

\bibitem{Sn20} V.P. Snaith:  Derived Langlands II: HyperHecke algebras, monocentric relatons and ${\mathcal M}_{cmc, \underline{\phi}}(G)$-admissibility; arXiv:3100675 [math.RT] 24 Mar 2020.

\bibitem{Sn20b} V.P. Snaith: Derived Langlands III: PSH algebras and their numerical invariants; arXiv:3223311 [math.RT] 12 June 2020.

\bibitem{newindnotes}  Derived Langlands IV: Notes on ${\mathcal M}_{c}(G)$-induced representations; arXiv:2008.06325v1 [math.RT] 14 Aug 2020.

\bibitem{Sn20c} V.P. Snaith: Derived Langlands V: The Hopflike properties of the hyperHecke algebra; preprint on University of Sheffield homepage (27 May 2020).


\bibitem{AVZ81} A.V. Zelevinsky: {\em Representations of finite classical groups - a Hopf algebra approach}; Lecture Notes in Math. \#869, Springer-Verlag (1981).

 \end{thebibliography}
\end{document}